\documentclass[12pt]{amsart}
\usepackage{amsmath}
\usepackage{amssymb}
\usepackage{tabularx}
\usepackage{enumerate}
\usepackage{graphicx}
\usepackage{texdraw}
\usepackage{color}

\usepackage{pgfplots}
\usepackage{graphics}
\usepackage{graphicx}
\usepackage{subfigure}
\usepackage{mathrsfs}

\usepackage{amsfonts,amssymb,amsmath}
\usepackage{epsfig}
\usepackage{bm}

\newtheorem{Theorem}{Theorem}[section]

\numberwithin{equation}{section}
\numberwithin{figure}{section}

\begin{document}
	
	\title[]{Prey-Predator Model on Graphs}
	
	\author[Y. Hu]{Yuanyang Hu$^1$}
	\author[C. Lei]{Chengxia Lei$^2$}
	
	\thanks{$^1$ School of Mathematics and Statistics,
		Henan University, Kaifeng, Henan 475004, P. R. China.}

	\thanks{$^2$ School of Mathematics and Statistics, Jiangsu Normal University, 
		Xuzhou, 221116, Jiangsu Province, China.}
	
	\thanks{{\bf Emails:} {\sf yuanyhu@mail.ustc.edu.cn} (Y. Hu).
		{\bf Emails:} {\sf leichengxia001@163.com} (C. Lei).
	}

	\begin{abstract}
		In this paper, we study the Lotka-Volterra prey-predator models consisting of two species on graphs under Neumann boundary condition and the condition that there is no boundary condition. Our results show that the existence of solutions and the global stability of the unique constant equilibrium solution of each parabolic system.	
	\end{abstract} 
	

	\maketitle
	
	\section{Introduction}	
	
	Let $G=G(V,E)$ be a finite connected weighted graph, where $V$ denotes the vertex set and $E$ denotes the edge set.
For each edge $xy \in E$, We suppose that its weight $w_{xy}>0$ and that $w_{xy}=w_{yx}$. Let $\mu: V\to \mathbb{R}^{+}$ be a finite measure. For any function $u: V \to \mathbb{R}$, the Laplacian of $u$ is defined by 
$$
\Delta_{V} u(x)=\frac{1}{\mu(x)} \sum_{y \sim x} w_{y x}(u(y)-u(x)),
$$
where $y \sim x$ means $xy \in E$.

For a subgraph $\Omega$ of a graph $G=G(V,E)$, the (vertex) boundary $\partial{\Omega}(\neq \emptyset)$ of $\Omega$ is the set of all vertices $z\in V$ not in $\Omega$ but adjacent to some vertex in $\Omega$, i.e., 
$$
\partial{\Omega}:=\{z\in V \backslash \Omega | z\sim y ~\text{for some}~y\in \Omega \}.
$$
Denote $\bar{\Omega}$ a graph whose vertices and edges are in $\Omega$ and vertices in $\partial{\Omega}$. For any function $u: \bar{\Omega} \to \mathbb{R}$, the Laplacian of $u$ on $\bar{\Omega}$ is defined by 
$$
\Delta_{\Omega} u(x)=\frac{1}{\mu(x)} \sum_{y\in \bar{\Omega}:y \sim x} w_{y x}(u(y)-u(x)).
$$
Let the (outward) normal derivative $\frac{\partial u}{\partial_{\Omega} n}$ at $z\in\partial \Omega$ by 
\begin{equation}\label{23}
\frac{\partial u }{\partial_{\Omega} n}(z):=\sum_{y\in {\Omega}} (u(z)-u(y))\frac{\omega_{zy}}{\mu(z)}.
\end{equation}

For given functions $u_0, v_0: V\to \mathbb{R}$ and $\tilde{u}_{0},\tilde{v}_{0}: \bar{\Omega}\to \mathbb{R}$, we study the following Lotka-Volterra prey-predator systems	
\begin{equation}\label{p}
\begin{cases}u_t-d_1 \Delta_{V} u=u\left(a_1-b_1 u-c_1 v\right), & x \in V,t>0,\\
 v_t-d_2 \Delta_{V} v=v\left(a_2+b_2 u-c_2 v\right), & x \in V, t>0, \\
  u(x, 0)=u_0(x)\ge 0,~ v(x, 0)=v_0(x)\ge 0, & x \in V,\end{cases}
\end{equation}
\begin{equation}\label{n}
\begin{cases}u_t-d_1 \Delta_{{\Omega}} u=u\left(a_1-b_1 u-c_1 v\right), & x \in \Omega, t>0, \\
 v_t-d_2 \Delta_{{\Omega}} v=v\left(a_2+b_2 u-c_2 v\right), & x \in \Omega, t>0, \\ 
 u(x, 0)=\tilde{u}_0(x)\ge 0,~  v(x, 0)=\tilde{v}_0(x)\ge 0, & x \in \bar{\Omega} ,\\
\frac{\partial u}{\partial_{\Omega} n}=0,~\frac{\partial v}{\partial_{\Omega} n}=0,& x \in \partial\Omega,
\end{cases}
\end{equation}
where the positive constants $d_i$ and $a_i$ ($i=1,2$) stand for the diffusion, growth rates of the populations $u$, $v$, respectively; $b_1>0$ and $c_2>0$ represent the rates of death  of the species $u$ and $v$, respectively; $c_1$ and $b_2$ are positive constants that measure the interaction of the populations $u$ and $v$, respectively.

	The system \eqref{n} can be regarded as the space-discretized version of the reaction diffusion system with Neumann boundary condition
	\begin{equation}\label{o}
		\begin{cases}u_t-d_1 \Delta u=u\left(a_1-b_1 u-c_1 v\right),&x\in Q,~t>0 \\ 
		v_t-d_2 \Delta v=v\left(a_2+b_2 u-c_2 v\right),&x\in Q,~t>0  ,\\
			u(x, 0)=\tilde{u}_0(x),~  v(x, 0)=\tilde{v}_0(x), &x \in \bar{Q} ,\\
			\frac{\partial u}{\partial n}=\frac{\partial v}{\partial n}=0,&x\in \partial Q.
		\end{cases}
	\end{equation} 
	where $Q\subset \mathbb{R}^{n}$ is a bounded domain, and $\Delta$ is the standard Laplacian operator, which describe the isotropic dispersal of species. The problem \eqref{o} has been studied by \cite{L} with $a_2$ replaced with a constant $e_2<0$. The authors in \cite{L} showed that all positive solutions of \eqref{o}, regardless of the initial value, converge to the constant equilibrium solution as $t\to \infty$. For more studies on diffusive Lotka-Volterra prey-predator systems of two species, one may refer to \cite{DS,DRS,LG,P,RZ} and the references therein. 
	
	In recent decade years, increasing efforts have been devoted to the investigation of anisotropy diffusion on graphs \cite{B,CB,GL,H,HL,LT}. For further related works, Gao et al. \cite{GL} use the method of establishing Lyapunov functions based on graph-theoretical approach to analyse the global stability for a predator-prey model among $n$ patches.
	Yokoi et al. \cite{Y} present two-patches system composed of prey and ambush predator. Nagatani \cite{N} deal with a metapopulation model on various graphs for a prey–predator system with different migration rates. 
	Slav\'{i}k \cite{S} consider a model of two competing species of Lotka-Volterra type with migration, where the spatial domain is an arbitrary network. They study the existence of stationary states and the asymptotic behavior of solutions.
	
	Here, one of our purpose of this paper is to discuss the boundeness of solutions to \eqref{p} and \eqref{n} by using upper and lower solutions method. Therefore, we have the following results:
	\begin{Theorem}\label{thm20221014}
		Let $G=(V,E)$ be a finite connected graph.	Suppose that $(u_0,v_0)\ge (0,0)$. Then \eqref{p} admits a unique solution $(u_1, u_2)$ that satisfies 
		\begin{equation}\label{b}
		(0,0) \leq(u_1, u_{2}) \leq\left(M_1, M_2\right)~\text{for}~x\in V~\text{and}~t\ge 0,
		\end{equation} 
		where  $$M_1=\max \{ \max\limits_{V} u_0(x), \frac{a_1}{b_1} \}\;\mbox{and}\;M_2=\max \left\{\frac{a_2+b_2 M_1}{c_2}, \max\limits _V v_0(x)\right\}.$$
	 Furthermore, if $u_0\not \equiv 0$ and $v_{0}\not \equiv 0$ on $V$, then $(u,v)$ is positive for  $x\in V$ and $t>0$.
	\end{Theorem}

	\begin{Theorem}\label{thm20221014-1}
		Let $\bar{\Omega}$ be a finite connected graph. Assume that $(\tilde{u}_0,\tilde{v}_0)\ge (0,0)$. Then \eqref{n} admits a unique solution $(u_3, u_4)$ that satisfies 
		\begin{equation}\label{b2}
		(0,0) \leq(u_3, u_{4}) \leq\left(M_3, M_4\right)~\text{for}~x\in \bar{\Omega}~\text{and}~t\ge 0,
		\end{equation} 
		where  $$M_3=\max \{ \max\limits_{\overline{\Omega}} \tilde{u}_0(x), \frac{a_1}{b_1} \}\;\mbox{and}\; M_4=\max \left\{\frac{a_2+b_2 M_3}{c_2}, \max\limits _{\overline{\Omega}} \tilde{v}_0(x)\right\}.$$
		 Moreover, if $\tilde{u}_0\not \equiv 0$ and $\tilde{v}_{0}\not \equiv 0$ on $\bar{\Omega}$, then $(u,v)$ is positive for $x\in \bar{\Omega}$ and $t>0$.
	\end{Theorem}
The equilibrium problems corresponding to \eqref{p} and \eqref{n} are the following elliptic systems:
\begin{equation}\label{p-1}
	\begin{cases}
		-d_1 \Delta_V u=u\left(a_1-b_1 u-c_1 v\right), & x \in V, \\ 
		-d_2 \Delta_V v=v\left(a_2+b_2 u-c_2 v\right), & x \in V, 
	\end{cases}
\end{equation}
and
\begin{equation}\label{n-1}
	\begin{cases}
		-d_1 \Delta_{\Omega} u=u\left(a_1-b_1 u-c_1 v\right), & x \in \Omega, \\ 
		-d_2 \Delta_{\Omega} v=v\left(a_2+b_2 u-c_2 v\right), & x \in \Omega, \\
		\frac{\partial u}{\partial_{\Omega} n}=0,~\frac{\partial v}{\partial_{\Omega} n}=0,& x \in \partial\Omega,
	\end{cases}
\end{equation}
respectively.
The biologically interesting case occurs when \begin{equation}\label{H1}
	\frac{a_1}{c_1}>\frac{a_2}{c_2}.
\end{equation} Clearly, under the hypothese \eqref{H1}, both \eqref{p} and \eqref{n} admit a unique positive constant steady-state $(e,g)$, where 
\begin{equation}\label{e}
e=\frac{a_1 c_2-c_1 a_2}{b_1 c_2+c_1 b_2}, \quad g=\frac{b_1 a_2+a_1 b_2}{b_1 c_2+c_1 b_2}.
\end{equation}	
The other goal of the present paper is to analyse the global stability of the equilibrium	$(e,g)$ by energy method. Now, we state the major results of the paper as following:
	\begin{Theorem}\label{11}
		Let $G=(V,E)$ be a finite connected graph and $(u_1,u_2)$ be the solution to \eqref{p}. If $\frac{a_1}{c_1}>\frac{a_2}{c_2}$, then 
		$$\lim\limits_{t\to \infty}(u_1,u_2)=(e,g)\;\mbox{ uniformly for}\; x\in V$$
provided that $u_0\ge, \not \equiv 0$ and $v_{0}\ge, \not \equiv 0$ on $V$, where $(e,g)$ is defined by \eqref{e}.
	\end{Theorem}
	\begin{Theorem}\label{22}
		Let $\bar{\Omega}$ be a finite connected graph and $(u_3,u_4)$ be the solution to \eqref{n}. If $\frac{a_1}{c_1}>\frac{a_2}{c_2}$, then  
		$$\lim\limits_{t\to \infty}(u_3,u_4)=(e,g)\;\mbox{uniformly for}\;x\in \bar{\Omega}$$
		provided that $\tilde{u}_0\ge, \not \equiv 0$ and $\tilde{v}_{0}\ge, \not \equiv 0$ on $\overline{\Omega}$, where $(e,g)$ is defined by \eqref{e}.
	\end{Theorem}
	The rest of the paper is organized as below. In section 2, we are devoted to  the proof of Theorems \ref{thm20221014} and \ref{thm20221014-1}. In  section 3, we give the proof of Theorems \ref{11} and \ref{22}.
	\section{The proof of Theorems \ref{thm20221014} and \ref{thm20221014-1}}
In this section, we are devoted to proving the existence and boundness of solutions of systems \eqref{p} and \eqref{n}. Now, we start with the proof of Theorem \ref{thm20221014}.

{\bf Proof of Theorem \ref{thm20221014}.} Suppose that $u_0\ge 0$ and $v_{0}\ge  0$ on $V$.
	Let $$M_1=\max \{ \max\limits_{V} u_0, \frac{a_1}{b_1} \}~\text{ and}~  M_2=\max \left\{\frac{a_2+b_2 M_1}{c_2}, \max\limits _V v_0(x)\right\}.$$ It is easy to check that $\left(\bar{u}_1, \bar{u}_2\right) \equiv\left(M_1, M_2\right)$ and $\left(\underline{u}_1, \underline{u}_2\right) \equiv(0,0)$ satisfy 
\begin{equation*}
	\begin{cases}	\begin{aligned}
			&\frac{\partial \bar{u}_1}{\partial t}-d_1 \Delta_V \bar{u}_1 \geq \bar{u}_1\left(a_1-b_1 \bar{u}_1-c_1 \underline{u}_2\right), x \in V, t>0, \\
			&\frac{\partial \bar{u}_2}{\partial t}-d_2 \Delta_V \bar{u}_2 \geq \bar{u}_2\left(a_2+b_2 \bar{u}_1-c_2 \bar{u}_2\right), x \in V, t>0, \\
			&\frac{\partial \underline{u}_1}{\partial t}-d_1 \Delta_V \underline{u}_1 \leq \underline{u}_1\left(a_1-b_1 \underline{u}_1-c_1 \bar{u}_2\right), x \in V, t>0, \\
			&\frac{\partial \underline{u}_2}{\partial t}-d_2 \Delta_V \underline{u}_2 \leq \underline{u}_2\left(a_2+b_2 \underline{u}_1-c_2 \underline{u}_2\right), x \in V, t>0,\\
			&\bar{u}_1(0) \geq u_0(x) \geq \underline{u}_1(0), \quad \bar{u}_2(0) \geq v_0(x) \geq \underline{u}_2(0),~x\in V.
	\end{aligned} 	\end{cases}
\end{equation*}
By \cite[Definition 3.4 ]{HL2}, we see that $(M_1,M_2)$ and $(0,0)$ is a pair of coupled upper and lower solutions of \eqref{p}. By virtue of \cite[Theorem 3.4]{HL2}, \eqref{p} admits a unique solution $(u_1, u_2)$ satisfying \begin{equation}\label{}
	(0,0) \leq(u_1, u_{2}) \leq\left(M_1, M_2\right)~\text{for}~x\in V~\text{and}~t\ge 0.
\end{equation}  
Furthermore, if $u_0\not \equiv 0$ and $v_{0}\not \equiv 0$ on $V$, then $u_1>0$ and $u_2>0$ for $x\in V$ and $t>0$, by \cite[Theorem 2.4]{HL2}.

{\bf Proof of Theorem \ref{thm20221014-1}.}	
Suppose that $\tilde{u}_0\ge 0$ and $\tilde{v}_{0}\ge  0$ on $\bar{\Omega}$.
Since $(\bar{u}_{3},\bar{u}_{4}) \equiv (M_3,M_4)$ and $(\underline{u}_{3},\underline{u}_{4})\equiv(0,0)$ satisfy 
the following system of inequalities
\begin{equation*}
	\begin{cases}	\begin{aligned}
			&\frac{\partial \bar{u}_3}{\partial t}-d_1 \Delta_{\Omega} \bar{u}_3 \geq \bar{u}_3\left(a_1-b_1 \bar{u}_3-c_1 \underline{u}_4\right), x \in \Omega, t>0, \\
			&\frac{\partial \bar{u}_4}{\partial t}-d_2 \Delta_\Omega \bar{u}_4 \geq \bar{u}_4 \left(a_2+b_2 \bar{u}_3-c_2 \bar{u}_4\right), x \in \Omega, t>0, \\
			&\frac{\partial \underline{u}_3}{\partial t}-d_1 \Delta_\Omega \underline{u}_3 \leq \underline{u}_3\left(a_1-b_1 \underline{u}_3-c_1 \bar{u}_4 \right), x \in \Omega, t>0, \\
			&\frac{\partial \underline{u}_4}{\partial t}-d_2 \Delta_\Omega \underline{u}_4 \leq \underline{u}_4\left(a_2+b_2 \underline{u}_3-c_2 \underline{u}_4\right), x \in \Omega, t>0,\\
			&\bar{u}_3(0) \geq \tilde{u}_0(x) \geq \underline{u}_3(0), \quad \bar{u}_4(0) \geq \tilde{v}_0(x) \geq \underline{u}_4(0),~x\in\bar{\Omega},
	\end{aligned} 	\end{cases}
\end{equation*}
By \cite[Definition 3.3]{HL2},  we can conclude that $(M_3,M_4)$ and $(0,0)$ is a pair of coupled upper and lower solutions of \eqref{n}. 
In view of \cite[Theorem 3.2]{HL2}, \eqref{n} admits a unique solution $(u_3, u_4)$ satisfying 
\begin{equation}\label{}
	(0,0) \leq(u_3, u_{4}) \leq\left(M_1, M_2\right)~\text{for}~x\in \bar{\Omega}~\text{and}~t\ge 0.
\end{equation} 
In addition, if we assume that $\tilde{u}_0\not \equiv 0$ and $\tilde{v}_{0}\not \equiv 0$ on $\bar{\Omega}$, then by Theorem 2.2 in \cite{HL2}, we see that $u_3>0$ and $u_4>0$ for $x\in \bar{\Omega}$ and $t>0$.

	\section{ The proof of Theorems \ref{11} and \ref{22} }
In this section, we prove the global stability of the unique positive constant steady-state $(e,g)$. For this purpose, we introduce some notations on graphs. For any functions $u$ and $v$ $: V\to \mathbb{R}$.  Define the gradient form of $u$ and $v$ as
	\begin{equation*}
	\Gamma(u, v)(x)=\frac{1}{2 \mu(x)} \sum_{y \sim x} w_{x y}(u(y)-u(x))(v(y)-v(x)).
	\end{equation*} Write $\Gamma(u)=\Gamma(u,u)$.
	Denote the length of the gradient of $u$ by
	\begin{equation*}
	|\nabla u|(x)=\sqrt{\Gamma(u)(x)}=\left(\frac{1}{2 \mu(x)} \sum_{y \sim x} w_{x y}(u(y)-u(x))^{2}\right)^{1 / 2}.
	\end{equation*}
	Define, for any function $
	u: V \rightarrow \mathbb{R}
	$, an integral of $u$ on $V$ by 
	$$\int \limits_{V} u d \mu=\sum\limits_{x \in V} \mu(x) u(x).$$	
For any functions $u,v$ on $V$, we use integration by parts to conclude that
	\begin{equation}\label{fb}
	\int_V u\Delta_{V} v d\mu=-\int_{V} \Gamma(u,v) d\mu.
	\end{equation}	
{\bf Proof of theorem \ref{11}.} For any solution $(u_1,u_2)$ of system \eqref{p}, we define 
		$$L(t):=\int\limits_{V} E(u_1(x,t),u_{2}(x,t)) d\mu$$ 
		where
		$$	E\left(u_1, u_2\right):= b_2 u_1-b_2 e \ln \left(\frac{u_1}{b_2 e}\right)+c_1 u_2-c_1 g \ln \left(\frac{u_2}{c_1 g}\right) .$$
	Clearly, \begin{equation}\label{c}
		E\left(u_1, u_2\right) \geq E(e, g)~\text{for~all}~ u_{1}>0~ \text{and}~ u_{2}>0.  
		\end{equation} 
		Then, we obtain \begin{equation*}
		\begin{aligned}
		\frac{dL}{d t} 
		=& \int_V\left[ b_2\left(1-\frac{e}{u_1}\right) \frac{\partial u_1}{\partial t} +c_1\left(1-\frac{g}{u_2}\right)  \frac{\partial u_2}{\partial t} \right]d\mu \\
		=& \int_V\bigg\{ b_2\left(1-\frac{e}{u_1}\right)\left[d_1 \Delta_V u_1+u_1\left(a_1-b_1 u_1-c_1 u_2\right)\right] + \\
		& c_1\left(1-\frac{g}{u_2}\right)\left[d_2 \Delta_V u_2+u_2\left(a_2+b_2 u_1-c_2 u_2\right)\right]\bigg\} d\mu.
		\end{aligned}
		\end{equation*}
		
		A direct calculation shows that  
		\begin{equation}\label{s1}
			\begin{aligned}
				& b_2\left(1-\frac{e}{u_1}\right) u_1\left(a_1-b_1 u_1-c_1 u_2\right)+c_1\left(1-\frac{g}{u_2}\right) u_2\left(a_2+b_2 u_1-c_2 u_2\right) \\
				=& b_2\left(u_1-e\right)\left[a_1-b_1\left(u_1-e\right)-c_1\left(u_2-g\right)-b_1 e-c_1 g\right] \\
				+& c_1\left(u_2-g\right)\left[a_2+b_2\left(u_1-e\right)-c_2\left(u_2-g\right)+b_2 e-c_2 g\right] \\
				=&-b_1 b_2\left(u_1-e\right)^2-c_1 c_2\left(u_2-g\right)^2.
			\end{aligned}
		\end{equation}
	From \eqref{b} and \eqref{fb}, we obtain that 
		\begin{equation*}
			\begin{aligned}
				\int\limits_V \left(1-\frac{e}{u_1}\right) \Delta_V u_1 d \mu &=\int\limits_V \Gamma\left(\frac{e}{u_1}, u_1\right) d \mu \\
				&=\int\limits_V \frac{1}{2 \mu (x)} \sum_{y \sim x}\left[\frac{e}{u_1(x)}-\frac{e}{u_1(y)}\right] w_{x y}\left[u_1(x)-u_1(y)\right] d \mu \\
				&=\int\limits_V \frac{-e}{2 \mu (x)} \sum_{y_{\sim} x} w_{x y} \frac{\left[u_1(x)-u_1(y)\right]^2}{u_1(x) u_1(y)} d \mu \\
				& \leq \int\limits_V-\frac{e}{M_1^2} \frac{1}{2 \mu (x)} \sum_{y \sim x} w_{xy} \left[u_1(x)-u_1(y)\right]^2 d \mu =-\frac{e}{M_1^2} \int\limits_{V} \left|\nabla u_1\right|^2 d \mu.
			\end{aligned}
		\end{equation*}
		Similarly, we have 
		\begin{equation*}
			\int\limits_V \left(1-\frac{g}{u_2}\right) \Delta_V u_2 d \mu \leq-\frac{g}{M_2^2} \int\limits_V \left|\nabla u_2\right|^2 d \mu \text {. }
		\end{equation*}
		Therefore, we deduce that 
		\begin{equation}\label{l}
			\begin{aligned}
				\frac{dL}{d t} & \leq-\frac{b_2 e d_1}{M_1^2} \int\limits_V \left|\nabla u_1\right|^2 d \mu-\frac{c_1 g d_2}{M_2^2} \int\limits_V \left|\nabla u_2\right|^2 d \mu\\
				&\;\;\;\;-\int\limits_V b_1 b_2\left(u_1-e\right)^2 d \mu-\int\limits_{V} c_1 c_2\left(u_2-g\right)^2 d \mu \\
				& \leq-\int\limits_V b_1 b_2\left(u_1-e\right)^2 d \mu-\int\limits_{V} c_1 c_2\left(u_2-g\right)^2 d \mu. 
			\end{aligned}
		\end{equation}

	Setting $F(t)=\int_V \left[u_1(x, t)-e\right]^2 d \mu$. It follows from \eqref{fb} that 
		\begin{equation*}
			\begin{aligned}
				\frac{d F}{d t} &=\int\limits_V 2\left(u_1(x, t)-e\right)\left[d_1 \Delta_V u_1+u_1\left(a_1-b_1 u_1-c_1 u_2\right)\right] d \mu \\
				&=2 \int\limits_V\left[-d_1 \Gamma\left(u_1, u_1\right)+\left(u_1-e\right) u_1\left(a_1-b_1 u_1-c_1 u_2\right)\right] d \mu .
			\end{aligned}
		\end{equation*}
		Thus, by \eqref{b}, we can find a constant $B=B(M_1,M_2,e,V)>0$ such that 
		\begin{equation}\label{4}
			\frac{dF}{dt}<B~\text{for}~t\ge 0.
		\end{equation}
		
		Now, we claim that
		\begin{equation}\label{F}
			\lim\limits_{t\to \infty} F(t)=0.
		\end{equation}
		Otherwise, we can find $\epsilon>0$ and a sequence $\{s_{n}\}_{n=1}^{\infty}$ satisfying $s_{n}+1<s_{n+1}$ and $s_{n}\to +\infty$ so that 
		\begin{equation*}\label{s}
			F(s_{n})>\epsilon~\text{for}~n=1,2,\cdots.
		\end{equation*} 
		Let $s=\min\{1,\frac{\epsilon}{2B}\}$, we use \eqref{4} to observe that 
		\begin{equation}\label{5}
			\frac{\varepsilon}{2}<F\left(s_n\right)-B s<F(t) \quad \text { for all } t \in\left[s_n-s, s_n\right], n=1,2, \cdots
		\end{equation}
		By virtue of \eqref{l} and \eqref{5}, we conclude that 
		$$
		\frac{dL}{d t} \leq \int_V-b_1 b_2\left(u_1-e\right)^2 d \mu=-b_1 b_2 F(t) \leq-b_1 b_2 \frac{\varepsilon}{2}
		$$
		for all $t \in\left[s_n-s, s_n\right], n=1,2, \cdots$.
		This implies that $L\left(s_n\right)-L\left(s_n-s\right) \leqslant-b_1 b_2 \frac{\varepsilon}{2} s$, which contradicts with \eqref{c}. Hence,  \eqref{F} holds. 
		
		Next, we assert that
		\begin{equation}\label{ug}
			u_1(x, t)-e \rightarrow 0 \text { as } t \rightarrow+\infty \text { for all } x \in V \text {. }
		\end{equation} 
		Otherwise, there exists $\epsilon>0$, $x_0\in V$ and a sequence $\{t_n\}_{n=1}^{\infty}$ such that $|u_1\left(x_0, t_n\right)-e|>\varepsilon, n=1,2, \cdots$. Thus we have 
		\begin{equation*}
			\begin{aligned}
				\int_V\left[u_1\left(x, t_n\right)-e\right]^2 d \mu &=\sum_{x \in V}\left[u_1\left(x, t_n\right)-e\right]^2 \mu(x) \\
				& \geq\left[u_1\left(x_0, t_n\right)-e\right]^2 \mu\left(x_0\right)>\varepsilon^{2} \mu\left(x_0\right)>0
			\end{aligned}
		\end{equation*}	
		which contradicts with \eqref{F} as $t\to \infty$.
		
		By a similar discussion, we see that 
		\begin{equation*}\label{9}
		\int_V\left[u_2(x, t)-g\right]^2 d \mu \rightarrow 0 \text { as } t \rightarrow+\infty \text {. }
		\end{equation*}
Then, we apply similar arguments to conclude that	
		\begin{equation*}
			u_2(x, t)-g \rightarrow 0 \text { as } t \rightarrow+\infty \text { for all } x \in V \text {. }
		\end{equation*}
		
		We now complete the proof.
		
{\bf Proof of Theorem \ref{22}}
		Define 
		$$\tilde{L}(t)=\int\limits_{\bar{\Omega}} E(u_{3}(\cdot,t),u_{4}(\cdot,t)) d\mu.$$ By \cite[Corollary 1.3 (ii)]{CB}, we know that
		$$	\int\limits_{\bar{\Omega}} \left(1-\frac{e}{u_j}\right) \Delta_{\Omega} u_j d \mu =\int\limits_{\bar{\Omega}} \Gamma\left(\frac{e}{u_j}, u_j\right) d \mu,\;j=3,4.$$ 
		From this and \eqref{s1}, we deduce that there exists constant $C_{1}>0$ such that
	$$	\frac{d\tilde{L}}{d t} \leq  -C_{1}\int\limits_{\bar{\Omega}}  \left(u_1-e\right)^2 d \mu-C_{1}\int\limits_{\bar{\Omega}}  \left(u_2-g\right)^2 d \mu.$$ 
		
 Define $$\tilde{F}(t):=\int\limits_{\bar{\Omega}} [u_{3}(\cdot,t)-e]^{2} d\mu.$$ Due to \cite[Corollary 1.3 (ii)]{CB}, we see that 
		\begin{equation*}
			\frac{d \tilde{F}}{d t} 
			=2 \int\limits_{\bar{\Omega}}\left[-d_1 \Gamma\left(u_3, u_3\right)+\left(u_3-e\right) u_3\left(a_1-b_1 u_3-c_1 u_4\right) \right]d \mu .
		\end{equation*}
		In view of \eqref{b2}, we see that $$\frac{d \tilde{F}}{d t}<\tilde{B},~t\ge 0$$ for some constant $\tilde{B}>0$. By a similar discussion as the proof of \eqref{F}, we see that $\lim\limits_{t\to \infty} \tilde{F}(t)=0$. Using a analogous argument as the proof of \eqref{ug}, we deduce that
		$$\lim\limits_{t\to \infty}u_{3}(x,t)=e\;\mbox{and}\;\lim\limits_{t\to \infty} u_{4}(x,t)=g \;\mbox{ for all}\; x\in \bar{\Omega}.$$
		
		We now complete the proof.

	\section*{Acknowledgments}
	
	Yuanyang Hu was partially supported by NSF of China (No. 12201184), National Natural Science Foundation of
	He' nan Province of China (No. 222300420416) and China Postdoctoral Science Foundation (No. 2022M711045).
Chengxia Lei was partially supported by NSF of China (No. 11801232, 11971454, 12271486), the NSF of
Jiangsu Province (No. BK20180999) and the Foundation of Jiangsu Normal University (No. 17XLR008).	
	
\end{document}